\documentclass[graybox]{svmult}

\usepackage{type1cm}        
\usepackage{makeidx}         
\usepackage{graphicx}        
\usepackage{multicol}        
\usepackage[bottom]{footmisc}

\usepackage{overpic}
\usepackage{newtxtext}       %
\usepackage{newtxmath}       
\usepackage{amscd}
\usepackage{wrapfig}
\usepackage[colorlinks = true,
            linkcolor = blue,
            urlcolor  = blue,
            citecolor = blue,
            anchorcolor = blue]{hyperref}

\makeindex             

\definecolor{airforceblue}{rgb}{0.36, 0.54, 0.66}
\definecolor{arsenic}{rgb}{0.23, 0.27, 0.29}
\definecolor{britishracinggreen}{rgb}{0.0, 0.26, 0.15}

\newcommand{\pp}[2]{\frac{\partial #1}{\partial #2}} 

\DeclareMathOperator{\diff}{d}

\newcommand{\LOh}    { \Lambda^0_h}

\newcommand{\LO}{      \Lambda^0                }  
\newcommand{\LI}{      \Lambda^1                }  
\newcommand{\LtwI}{    \widetilde\Lambda^1      } 
\newcommand{\LtwO}{    \widetilde\Lambda^0      } 
\newcommand{\LstkI}{    \Lambda^{k+1}      } 
\newcommand{\Lstk}{     \Lambda^{k}      } 
\newcommand{\LtwkI}{    \widetilde\Lambda^{k+1}      } 
\newcommand{\Ltwk}{    \widetilde\Lambda^{k}      } 
\newcommand{\hItw}{\widetilde h^{(1)}} 
\newcommand{\hOst}{h^{(0)}} 
\newcommand{\uIst}{u^{(1)}}
\newcommand{\uOtw}{\widetilde u^{(0)}}  
\newcommand{\vOst}{ v^{(0)}}  
\newcommand{\vItw}{\widetilde v^{(1)}}  
\newcommand{\qOtw}{\widetilde q^{(0)}}   
\newcommand{\FuOtw}{\widetilde F_u^{(0)}} 
\newcommand{\FvOst}{F_v^{(0)}} 
\newcommand{\uIhst}{u_h^{(1)}}
\newcommand{\vIhtw}{\widetilde v_h^{(1)}}
\newcommand{\wIhst}{ \chi_h^{(1)} }
\newcommand{\wIhtw}{ \widetilde \chi_h^{(1)} }
\newcommand{\uOhtw}{\widetilde u_h^{(0)}}  
\newcommand{\vOhst}{v_h^{(0)}}  
\newcommand{\hIhtw}{\widetilde h_h^{(1)}} 
\newcommand{\hOhst}{h_h^{(0)}}
\newcommand{\phiOhtw}{\widetilde \phi_h^{(0)}}
\newcommand{\chiIhtw}{\widetilde \chi_h^{(1)}} 
\newcommand{\BOhst}{B_h^{(0)}}
\newcommand{\BOst}{B^{(0)}}
\newcommand{\qOhtw}{\widetilde q_h^{(0)}}   
\newcommand{\FuOhtw}{\widetilde F_{u_h}^{(0)}} 
\newcommand{\FvOhst}{F_{v_h}^{(0)}} 
\newcommand{\Hh}{\mathcal{H}}
\newcommand{\Fh}{\mathcal{F}}
\newcommand{\Ch}{\mathcal{C}}
\newcommand{\Gh}{\mathcal{G}}

\newcommand{\dede}[2]{\frac{\delta #1}{\delta #2}}
\newcommand{\tstar}{\tilde \star_h }
\newcommand{\ttstar}{\tilde \star }
\newcommand{\indxVOtest} {  {   {l'}       } } 
\newcommand{\indxVOtrial}{  {   {l}      } } 
\newcommand{\indxVItest} {  {   {m'}       } } 
\newcommand{\indxVItrial}{  {   {m}      } } 
\newcommand{\fctVOtest}{   \phi_{l'}(x)      } 
\newcommand{\fctVOtrial}{   \phi_{l}(x)      }
\newcommand{\fctVItest}{   \chi_{m'}(x)      } 
\newcommand{\fctVItrial}{   \chi_{m}(x)      }
\newcommand{\hVOtrial}{   h_{l}      } 
\newcommand{\uVItrial}{   u_{m}      } 
 
\newcommand{\htwVItrial}{   \widetilde{h}_{m}      } 
\newcommand{\LIh}{     \Lambda^1_h              }
\newcommand{\GPIu}{ {\rm GP1}_{u} }  
\newcommand{\GPOu}{ {\rm GP0}_{u} }  
\newcommand{\GPIh}{ {\rm GP1}_{h} }  
\newcommand{\GPOh}{ {\rm GP0}_{h}  }  
\newcommand{\LtwIh}{   \widetilde\Lambda^1_h    }
\newcommand{\LtwOh}{   \widetilde\Lambda^0_h    }
\newcommand{\Mnn}{ {\bf M}^{nn} } 
\newcommand{\Mee}{ {\bf M}^{ee} } 
\newcommand{\Den}{ {\bf D}^{en} } 
\newcommand{\Dne}{ {\bf D}^{ne} }
\newcommand{\Mne}{ {\bf M}^{ne}  } 
\newcommand{\Men}{ {\bf M}^{en}  }

\newcommand{\prjcmxne}{ {\bf P}^{ne}  }  
\newcommand{\Pne}{ {\bf P}^{ne}  }  
  
\newcommand{\metrivcte}{  {\bf \Delta x}_{e}  }
\newcommand{\qnvct}{ {\bf \tilde q}^0_{n}  }
\newcommand{\uOn}{ {\bf \tilde u}^0_{n}  } 
\newcommand{\vOn}{ {\bf   v}^0_{n}  } 
\newcommand{\BOn}{ {\bf B}^0_{n} } 
\newcommand{\hOn}{ {\bf h}^0_{n}  } 
\newcommand{\uIe}{ {\bf u}^1_{e}  }
\newcommand{\vIe}{ {\bf   \tilde v}^1_{e}  }
\newcommand{\hIe}{ {\bf \tilde h}^1_{e}  }
\newcommand{\FuOntw}{ {\bf   \tilde F}^{u}_n  }
\newcommand{\FvOnst}{ {\bf   F}^{v}_n  }
\newcommand{\AVGu}{ {\rm AVG}_{u} }  
\newcommand{\AVGh}{ {\rm AVG}_{h} }

\usepackage{lineno}

\begin{document}

 \title*{A structure-preserving approximation of the discrete split rotating shallow water equations}
\titlerunning{A structure-preserving approximation of the discrete split RSW equations}
\author{Werner Bauer, J\"orn Behrens, Colin J. Cotter}
\authorrunning{W. Bauer et al.}
\institute{Werner Bauer \at Inria Rennes, France, and Imperial College London, UK, \href{mailto:werner.bauer.email@gmail.com}{werner.bauer.email@gmail.com}
\and J\"orn Behrens \at CEN/Depart. of Mathematics, Universit\"at Hamburg, DE, \href{mailto:joern.behrens@uni-hamburg.de}{joern.behrens@uni-hamburg.de}
\and Colin J. Cotter \at Imperial College London, London SW7 2AZ, UK, \href{mailto:colin.cotter@imperial.ac.uk}{colin.cotter@imperial.ac.uk}
}
%
%
\maketitle

\abstract{
We introduce an efficient split finite element (FE) discretization of 
a y-independent (slice) model of the rotating shallow water equations.
The study of this slice model provides insight towards developing schemes 
for the full 2D case. Using the split Hamiltonian FE framework \cite{BaBeCo2019} 
(Bauer, Behrens and Cotter, 2019), we result in structure-preserving discretizations 
that are split into topological prognostic and metric-dependent closure equations. 
This splitting also accounts for the schemes' properties: the Poisson bracket 
is responsible for conserving energy (Hamiltonian) as well as mass, potential 
vorticity and enstrophy (Casimirs), independently from the realizations of the 
metric closure equations. The latter, in turn, determine accuracy, stability, 
convergence and discrete dispersion properties. We exploit this splitting to 
introduce structure-preserving approximations of the mass matrices in the metric 
equations avoiding to solve linear systems. We obtain a fully structure-preserving 
scheme with increased efficiency by a factor of two.
}

\vspace{-1em}
\section{Introduction}
\label{sec_intro}
\vspace{-1em}

The notion of structure-preserving schemes describes discretizations that preserve 
important structures of the corresponding continuous equations: e.g.
(i) the conservation of invariants such as energy, mass, vorticity and enstrophy in 
the case of the rotating shallow water (RSW) equations, 
(ii) the preservation of geometric structures such as {\sf div curl} = {\sf curl grad} = 0 
or the Helmholtz decomposition of vector fields, and 
(iii) the conservation of large scale structures such as geostrophic or hydrostatic 
balances \cite{Staniforth2012}. 
Their conservation is important to avoid, for instance, biases in the statistical behavior 
of numerical solutions \cite{Franck2007} or to get models that correctly transfer energy 
and enstrophy between scales \cite{NataleCotter2017b}.

The construction of such schemes is an active area of research and 
various approaches to develop structure-preserving discretizations exist: 
e.g. variational discretizations \cite{BaGB2018RSW,BaGB2017AnPI,PaMuToKaMaDe2010}
or compatible FE methods \cite{CotterThuburn2012,McRaeCotter2014}. 
In particular FE methods are a very general, widely applicable approach 
allowing for flexible use of meshes and higher order approximations. 
When combined with Hamiltonian formulations, they allow
for stable discretizations of the RSW equations that conserve energy and enstrophy 
\cite{BauerCotter2018,McRaeCotter2014}. 
However, they usually apply integration by parts to address 
the regularity properties of the FE spaces in use, which introduces additional 
errors and non-local differential operators. Moreover, FE discretizations 
usually involve mass matrices which are expensive to solve, while
approximations of the mass matrices have to be designed carefully in order to preserve structure.

To address these disadvantages, we introduced in \cite{Bauer2017_1D,BaBeCo2019} 
the split Hamiltonian FE method based on the split equations of GFD \cite{Bauer2016}, 
in which pairs of FE spaces are used such that integration by parts is avoided,
and we derived structure-preserving discretizations of a y-independent RSW slice-model.
Our method shares some basic ideas with {\em mimetic discretizations} 
(e.g. \cite{Beirao2017,Bochev2006,CotterThuburn2012,PalhaGerritsma2014})
in which PDEs are written in differential forms, but stresses a distinction between
topological and metric parts and the use of a proper FE space for each variable.

 Here, we address the disadvantage of FE methods arising from mass matrices. 
 In the framework of split FEM \cite{Bauer2017_1D,BaBeCo2019}, we introduce approximations
 of the mass matrices in the metric equations resulting in a structure-preserving
 discretization of the split RSW slice-model that is more efficient that the original 
 schemes introduced in \cite{BaBeCo2019}. To this end, we recall in Sect.~\ref{sec_splitHam}
 the split Hamiltonian framework and the split RSW slice-model, and we introduce the 
 approximation of the metric equations. In Sect.~\ref{sec_numerics}, we present 
 numerical results and in Sect.~\ref{sec_conclusion} we draw conclusions.

\vspace{-2em}
\section{Split Hamiltonian FE discretization of the RSW slice-model}
\label{sec_splitHam}
\vspace{-1em}

On the example of a y-independent slice model of the RSW equations,
we recall the split Hamiltonian FE method of \cite{BaBeCo2019}. 
For pairs of height fields (straight 0-form $\hOst$ and twisted 1-form $\hItw$),
of velocity fields in $x$-direction (twisted 0-form $\uOtw$ and straight 1-form $\uIst$),
and of velocity fields in outer slice direction (straight 0-form $\vOst$ and twisted 1-form $\vItw$),
this RSW slice-model reads
 \begin{equation}\label{slice-model}
 \begin{split}
  \pp{\uIst}{t} - \tilde\star \qOtw  \FvOst    + \diff \BOst  = 0, & \quad   
  \pp{\vItw}{t} + \tilde\star \qOtw  \FuOtw     = 0,\quad
       \pp{\hItw}{t}    +  {\diff} \FuOtw    = 0,  \\ 
   \uOtw  = \widetilde\star \uIst  , & \qquad  
   \vOst  = \widetilde\star \vItw, \qquad  
   \hItw  = \widetilde\star \hOst  , 
 \end{split}
 \end{equation}
 in which $\FuOtw := \hOst \uOtw$ and $\FvOst := \hOst \vOst$
 are mass fluxes, $\BOst := g \hOst + \frac{1}{2}(\uOtw)^2 + \frac{1}{2}(\vOst )^2$
 is the Bernoulli function with gravitational constant $g$, and 
 $\qOtw   \hItw = {\diff} \vOst + f dx $ with Coriolis 
 parameter $f$ defines implicitly the potential vorticity (PV) $\qOtw$. 
 All variables are functions of $x$ and $t$: for instance,  
 $u(x,t)$ is the coefficient function of the 1-form $\uIst = u(x,t) dx$.

  \begin{wrapfigure}{r}{.5\textwidth}
  \vspace{-1em}
  \begin{equation}\notag
   \begin{CD}  
     \hOst ,  \vOst \  \in \  \LO    @>\diff >>  \LI    \  \ni \  \uIst   \\
              \hspace*{+1.6cm}    @V  \tilde \star  \, VV            \hspace*{-1.1cm}               @VV \, \tilde \star V          \\
     \hItw , \vItw \  \in \  \LtwI   @<\diff <<  \LtwO  \  \ni \  \uOtw   \\ 
    \end{CD}
   \end{equation}
   \caption{Relation between operators and spaces.}
   \label{CD_1}
   \vspace{-1em}
 \end{wrapfigure}
 The pairs of variables are connected via the 
 twisted Hodge-star operator $\widetilde \star: \Lambda^k \rightarrow \widetilde \Lambda^{(1 - k)}$ 
 that maps from straight $k$-forms to twisted $(1-k)$-forms (or vice versa) with $k=0,1$ in one dimension (1D). 
 The index $^{(k)}$ denotes the degree, and $\Lambda^k,\widetilde \Lambda^k$ the space of all $k$-forms. 
 Note that straight forms do not change their signs when the orientation 
 of the manifold changes in contrast to twisted forms. 
 The exterior derivative $\diff$ is a mapping $\diff: \Lstk \rightarrow \LstkI$ 
 (or $\diff: \Ltwk \rightarrow \LtwkI$).
 Here in 1D, it is simply the total derivative of a smooth function $g^{(0)} \in \LO, \diff g^{(0)} = \partial_x g(x) dx \in \LI$ 
 (see \cite{Bauer2016} for full details).
 Diagram~\eqref{CD_1} illustrates the relations 
 between the operators and spaces.

 \smallskip\noindent\textbf{Galerkin discretization.}
  To substitute FE for continuous spaces, we 
  consider $\LOh, \LtwOh = CG_p$ and $\LIh, \LtwIh = DG_{p-1}$ with polynomial order $p$. 
  The discrete Hodge star operators $\tstar^0,\tstar^1$ map between straight and twisted spaces 
  and may be non-invertible. 
  The split FE discretization of Eqns.~\eqref{slice-model} 
   seeks solutions $\uIhst,\vIhtw, \hIhtw \in (\LIh (L),\LtwIh (L),\LtwIh (L))$ 
   of the topological equations (as trivial projections)
   \begin{align}
    \langle \wIhst, \pp{}{t} \uIhst \rangle - \langle \wIhst, \ttstar \qOhtw \FvOhst \rangle   
    + \langle \wIhst, \diff \BOhst  \rangle = 0, 
    & \qquad \forall \wIhst \in \LIh, \label{equ_weak_momtu} \\
    \langle \wIhtw, \pp{\vIhtw}{t}  \rangle + \langle \wIhtw, \ttstar \qOhtw \FuOhtw \rangle  = 0, 
    & \qquad \forall \wIhtw \in \LtwIh, \label{equ_weak_momtv} \\
    \langle \chiIhtw , \pp{\hIhtw}{t} \rangle   + \langle \chiIhtw , \diff \FuOhtw  \rangle = 0, 
    & \qquad \forall  \chiIhtw \in \LtwIh, \label{equ_weak_cont}\\
  \langle \ttstar \phiOhtw ,  \qOhtw \hIhtw \rangle - \langle \ttstar\phiOhtw ,  \diff \vOhst \rangle 
  -\langle \ttstar \phiOhtw ,  f dx\rangle =0, 
    & \qquad \forall \phiOhtw \in \LtwOh, \label{pv}
    \end{align}
  subject to the metric closure equations (as non-trivial Galerkin projections (GP))
  \begin{equation}\label{equ_varsplit_metri_disc} 
  \begin{split}
  \uOtw  = \widetilde\star \uIst \approx {\rm GP}(1-i)_u: \langle \uOhtw (x) , \hat\tau^i (x) \rangle & =  \langle  \uIhst (x),  \hat\tau^i (x)\rangle  \, , 
         \        \forall \hat\tau^i  \in \hat\Lambda_h^{i},    \\  
   \vOst  = \widetilde\star \vItw  \approx  {\rm GP}(1-i)_v:  \langle \vOhst (x)  , \hat\tau^i (x)  \rangle & =  \langle  \vIhtw (x),  \hat\tau^i (x) \rangle  \, , 
         \        \forall \hat\tau^i  \in \hat\Lambda_h^{i},    \\ 
     \hItw  = \widetilde\star \hOst  \approx  {\rm GP}(1-i)_h:  \langle \hOhst (x), \hat\tau^{j} (x) \rangle & = \langle  \hIhtw (x), \hat\tau^{j}(x)  \rangle \, ,  
          \         \forall \hat\tau^{j}  \in \hat\Lambda_h^{j},  
  \end{split}  
  \end{equation} 
  for $i,j = 0,1$. 
  $\langle \cdot ,\cdot \rangle := \int_L \cdot \wedge \ttstar \cdot$ is the 
  $L_2$ inner product on the domain $L$. 
  We distinguish between continuous and discrete Hodge star operators $\ttstar$ 
  and $\tstar$, respectively.  
  $\ttstar$ is used in $\langle, \rangle$ such 
  that $k$-forms of the same degree are multiplied, which 
  translates into a standard inner product for coefficient functions: e.g. 
  $\langle \ttstar \uIhst, \uOhtw \rangle =\langle \uIhst(x), \uOhtw(x)\rangle $.
  $\tstar$ is realized as in Eqns.~\eqref{equ_varsplit_metri_disc} 
  via non-trivial GP between 0- and 1-forms, cf. \cite{Bauer2017_1D}.

  \vspace{-2em}
  \subsection{Continuous and discrete split Hamiltonian RSW slice-model}
  \label{sec_splitham}   
   \vspace{-1em}

  Both the continuous split RSW slice-model of Eqns.~\eqref{slice-model} and the corresponding weak 
  (discrete) form of \eqref{equ_weak_momtu}--\eqref{equ_varsplit_metri_disc} can equivalently be written 
  in Hamiltonian form, as shown in \cite{BaBeCo2019}. 
  Considering the discrete version, the Hamiltonian with metric equations reads
\begin{equation}\label{disc_H}
\begin{split}
       \Hh[\uIhst,\vIhtw, \hIhtw] & =  \frac{1}{2} \langle   \uIhst , {\ttstar}\hOhst \uOhtw \rangle
				   + \frac{1}{2} \langle\vIhtw ,  {\ttstar}\hOhst \vOhst \rangle 
				   +  \langle  \hIhtw , {\ttstar} g \hOhst \rangle \\
				   &  \uOhtw =     \tilde \star_h \uIhst , \ \vOhst =  \tilde \star_h \vIhtw, \ \hOhst =   \tilde \star_h \hIhtw 
				   \  \text{(metric eqns.)}
\end{split}
\end{equation}
 while the almost \emph{Poisson bracket} $\{ , \}$ is defined as
 {\fontsize{8pt}{8pt}\selectfont
 \begin{equation}\label{eqn_bracket_split}
 \{ \Fh, \Gh \} := - \langle \dede{\Fh}{\hIhtw}, \diff \ttstar \dede{\Gh}{\uIhst} \rangle 
                  - \langle \dede{\Fh}{\uIhst}, \diff \ttstar \dede{\Gh}{\hIhtw} \rangle 
                  + \langle \dede{\Fh}{\uIhst}, \ttstar  {\qOhtw} \ttstar \dede{\Gh}{\vIhtw}\rangle 
                  - \langle \dede{\Fh}{\vIhtw}, \ttstar  {\qOhtw} \ttstar \dede{\Gh}{\uIhst}\rangle 
 \end{equation}
 }
 
 \vspace{-1em}
 \noindent with \emph{PV} defined by:
 $\langle \ttstar\phiOhtw ,   \qOhtw \hIhtw \rangle + \langle \diff \phiOhtw , \vIhtw \rangle 
  -\langle \ttstar\phiOhtw ,  f dx\rangle =0, 
     \ \forall \phiOhtw \in \LtwOh $.
 Then, the {dynamics} for any functional $\Fh[\uIhst,\vIhtw,\hIhtw]$ is given by 
$  \frac{d} {dt}\Fh = \{ \Fh, \Hh \}$.

\smallskip\noindent\textbf{Splitting of schemes properties.}
 The split Hamiltonian FE method results in a family of schemes in which 
 the schemes' properties split into \emph{topological and metric dependent 
 ones}, cf. \cite{BaBeCo2019}.
 
 The \textbf{topological properties} hold for all double FE pairs that fulfill commuting diagram~\eqref{CD_1}.
 In particular,  \vspace{-0.5em}
  \begin{itemize}
   \item the total energy is conserved, because $\frac{d}{dt} \Hh  = \{ \Hh, \Hh \}   = 0$\ 
   which follows from the antisymmetry of \eqref{eqn_bracket_split};
   \item the Casimirs $\Ch=M,PV,PE$ are conserved as 
   $\frac{d}{dt} \Ch  = \{ \Ch, \Hh \}   = 0$ for $\{\Ch, \Gh\} = 0 \forall \Gh$
   with $\Ch = \langle  \hIhtw ,{\ttstar} F(\qOhtw) \rangle $ 
   for $F=1(M)$, $F= \qOhtw\tilde 1(PV)$, $F=(\qOhtw)^2(PE)$; 
   
   \item \{, \} is independent of $\tstar$, hence $\Hh,\Ch$ are 
   conserved for any metric equation.
  \end{itemize}

  The \textbf{metric properties} are associated to a certain choice of FE spaces.
  In particular, this choice determines\vspace{-0.5em}
    \begin{itemize}
     \item the dispersion relation which usually depends on $\Delta x$ between DoFs,
     \item the stability, because the inf-sub condition depends on the norm, and 
     \item convergence and accuracy, which both are measured with respect to norms.
    \end{itemize}
    
  \vspace{-2em} 
  \subsection{Family of structure-preserving split RSW schemes}
  \vspace{-1em}

 Besides the splitting into topological and metric properties, another remarkable 
 feature of the split FE framework is that one choice
 of compatible FE pairs leads to a family of split schemes, cf. \cite{Bauer2017_1D,BaBeCo2019}. 
 In the following, we consider for $p=1$ the piecewise linear space $\LOh, \LtwOh = CG_p = {\rm P}1$ 
 with basis $\{ \fctVOtrial \}_{\indxVOtrial = 1}^{N}$ and the piecewise constant space 
 $\LIh, \LtwIh = DG_{p-1} = {\rm P}0$ with basis $\{ \fctVItrial \}_{\indxVItrial = 1}^{N}$.
 Being in a 1D domain with periodic boundary, both have 
 $N$ independent DoFs. We approximate 0-forms in P0 and 1-forms in P1, 
 e.g.  $\uIhst(x,t) = \sum_{\indxVItrial = 1}^{N} \uVItrial(t) \fctVItrial$.
 The split framework \cite{BaBeCo2019} leads to \emph{one set of discrete topological equations}
 for Eqn.~\eqref{equ_weak_momtu}--\eqref{pv}, and 
 \emph{four combinations of discrete metric equations} for \eqref{equ_varsplit_metri_disc}:
 \begin{small}
   \begin{eqnarray}\notag
 \hspace{0em}\text{\textbf{topol. momentum eqn.:}} \hspace{1em}  
   \pp{}{t} \uIe + \Den \BOn   - \widetilde{\Men} (\qnvct  \circ \FvOnst )  = 0 , \ 
        \pp{}{t} \vIe  + \widetilde{\Men} (\qnvct \circ \FuOntw  )  = 0 ,\notag\\
        \hspace{0em}\text{\textbf{topol. continuity eqn.:}} \hspace{20.7em}   \pp{}{t} \hIe + \Den \FuOntw = 0 ,\notag 
    \end{eqnarray}
    \begin{equation}\label{equ_metric}
   \hspace{-3.5em}\text{\textbf{metric closure eqn.:}}\hspace{1.5em}  
    \begin{CD}  
     \hOn     \quad \in \quad \LOh  {\subset P1}       @>\Den >>           \LIh {\subset P0}  \quad \ni \quad  \uIe ,
     (\vIe \in \LtwIh)\\
              \hspace*{+1.2cm}    @V  {{\GPIh:} \ \Mnn \hOn          = \prjcmxne        \hIe   }  VV     \hspace*{-3.2cm}   @VV   
                                     {{\GPIu:} \ \Mnn  \uOn          = \prjcmxne        \uIe \ \& \ \Mnn  \vOn          = \prjcmxne        \vIe }  V    \\
              \hspace*{+1.2cm}    @V  {{\GPOh:} \ \Men \hOn          = \hIe  \, \quad          }  VV     \hspace*{-3.2cm}   @VV   
                                     {{\GPOu:} \ \Men  \uOn          = \uIe   \ \ \, \& \ \, \Men  \vOn          = \vIe                 }  V     \\
    \hIe     \quad \in \quad \LtwIh  {\subset P0}      @<\Den  <<  \LtwOh {\subset P1} \quad \ni \quad \uOn , 
     (\vOn \in \LOh). \\ 
    \end{CD}
    \end{equation}
 \end{small}

\noindent
Here, we used the following $(N\times N)$ matrices 
with index $n$ for nodes and $e$ for elements: 
(i) mass matrices $\Mnn$, $\Mee$, $\Men$, with metric-dependent coefficients 
$ (\Mnn)_{\indxVOtrial\indxVOtest} = \int_L \fctVOtrial \fctVOtest  dx$,
$(\Mee)_{\indxVItrial\indxVItest} = \int_L \fctVItrial \fctVItest  dx $,
$(\Men)_{\indxVOtrial\indxVItest} =  \int_L \fctVOtrial  \fctVItest dx$
(with $\widetilde{\Men} = \{\Men \text{in Or} ,-\Men \text{in -Or}\}$ for orientation Or of $L$ 
and $\Men = (\Mne)^T$ with $T$ for the transposed matrix); 
and (ii) the stiffness matrix $\Den$ with metric-independent coefficient
$(\Den)_{\indxVOtrial\indxVItest} = \int_L  \frac{d \fctVOtrial  }{dx}\fctVItest  dx$
(with $\Den = (\Dne)^T$).
 We separate $\Mne = \Pne \,(\metrivcte)^T$ into a metric-dependent $\metrivcte$ 
 and a metric-free part $\prjcmxne$, the latter is an averaging operator 
 from $e$ to $n$ values. 
  
 Considering the discrete variables, $\uIe = \Mee {\bf u}_e$ is a discrete 1-form associated to the vector array 
 ${\bf u}_e = \{u_m(t)| m =1, ...N\}$, analogously for $\vIe$. 
 Similarly, $ \hIe  = \Mee {\bf \tilde h_e}$ is a discrete 1-from 
 with ${\bf \tilde h}_e = \{\htwVItrial (t)| m =1, ...N\}$. 
 Discrete 0-forms read, e.g. $\hOn = \{\hVOtrial (t)| l =1, ...N\}$.
 The discrete mass fluxes are 
 $\FuOntw = \hOn \circ \uOn \quad \text{and} \quad \FvOnst = \hOn \circ \vOn$
 and the discrete Bernoulli function reads
 $  \BOn = \frac{1}{2} \uOn \circ \uOn +  \frac{1}{2} \vOn \circ \vOn  + g \hOn$.
 Finally, $\GPIu,\GPOu,\GPIh,\GPOh$ are the nonlinear GPs of \eqref{equ_varsplit_metri_disc} 
 for P1 and P0 test functions.

   \vspace{-2em}
  \subsection{A structure-preserving approximation of split RSW schemes}
   \vspace{-1em}
 \begin{wrapfigure}{r}{.4\textwidth} 
 \vspace{-1em}
  \centering
  \includegraphics[scale=.78]{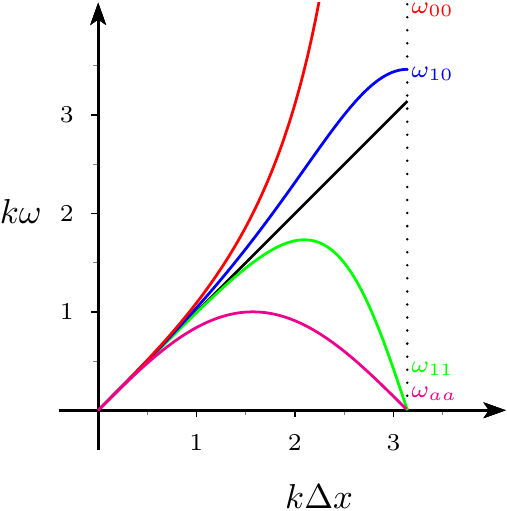}
 \caption{Dispersion relations: 
  analytic (black) for $c=\sqrt{gH}=1$, 
  $\omega_{00}$ (green) for $\GPIu$-$\GPIh$, 
  $\omega_{01}$ (blue) for $\GPIu$-$\GPOh$ and $\GPOu$-$\GPIh$,
  $\omega_{11}$ (red) for $\GPOu$-$\GPOh$ (cf. \cite{BaBeCo2019}), and 
    $\omega_{aa}$ (magenta) for ${\rm AVG}_u$-${\rm AVG}_h$.}
 \label{fig_disprel}
 \vspace{-1.5em}
  \end{wrapfigure}
Here we introduce a new, computationally more efficient split RSW scheme 
compared to those of \cite{BaBeCo2019}. 
We exploit the splitting of the topological and metric properties within 
the split FE framework to introduce structure-preserving approximations
of the mass matrices used in the metric equations. 
Instead of using the full nontrivial Galerkin projections 
$\GPIh,\GPOh$ for height or $\GPIu,\GPOu$ for velocity $u,v$, 
we use the \textbf{averaged} versions: 
$$
 {\rm AVG}_h: \hOn = \mathbf{A}^{en}   \hIe \, , \quad {\rm AVG}_u:\uOn = \mathbf{A}^{en}   \uIe \, ,
$$
with averaging operator $\mathbf{A}^{en}= (\Pne)^T$
and denote the resulting scheme with $\AVGu$-$\,\AVGh$. 
Rather then solving linear systems in \eqref{equ_metric}, 
we obtain values for $\hOn,\uOn,\vOn$ simply by averaging. 
This is computationally more efficient. 
In fact, already for this 1D problem we achieve a speedup by 
a factor of 2 (wall clock time) compared to the full GPs.

As stated in Sect.~\ref{sec_splitham}, such modification 
does not impact on the structure-preserving properties
but will change the metric-dependent ones instead. 
Before we confirm this in Sect.~\ref{sec_numerics} numerically,
we first determine analytically the discrete dispersion relation 
related to this approximation. 
A similar calculation as done in \cite{Bauer2017_1D} leads to the 
following discrete dispersion relation: 

\vspace{-0.8em}
$$  c_d =  \frac{\omega_{aa}}{k} = \pm \sqrt{gH}\frac{1}{k \Delta x} \sin (k \Delta x)$$
\vspace{-1.2em}

\noindent
with angular frequency $\omega_{aa} = \omega_{aa}(k)$ and discrete wave speed 
$c_d \rightarrow c = \sqrt{gH}$ (with mean height $H$)
in case $k \rightarrow 0$ and with a spurious mode (second zero root) at shortest 
wave length $k = \frac{\pi}{\Delta x}$. 
As shown in Fig.~\ref{fig_disprel} (with results relative to the nondimensional wave speed $c = \sqrt{gH}=1$), 
this is similar to the dispersion relation of the $\GPIu$-$\GPIh$ scheme in the sense 
that both have a spurious mode at $k = \frac{\pi}{\Delta x}$, cf. \cite{Bauer2017_1D}. 
For completeness, we added the dispersion relation for the other possible 
realizations of the metric equations \eqref{equ_metric} as introduced 
in \cite{Bauer2017_1D,BaBeCo2019}.

\vspace{-2em} 
\section{Numerical results}
\label{sec_numerics}
\vspace{-1em}

We study the structure-preserving properties, as well as convergence, stability and 
dispersion relation for the averaged split scheme $\AVGu$-$\,\AVGh$ and compare it with the split schemes 
of \cite{BaBeCo2019}. We use test cases (TC) in the quasi-geostrophic regime 
such that effects of both gravity waves and compressibility are important.

\begin{figure}[h!] \centering \vspace{-1em}
  \begin{tabular}{ccc} 
    \hspace{0cm}{\includegraphics[scale=0.35]{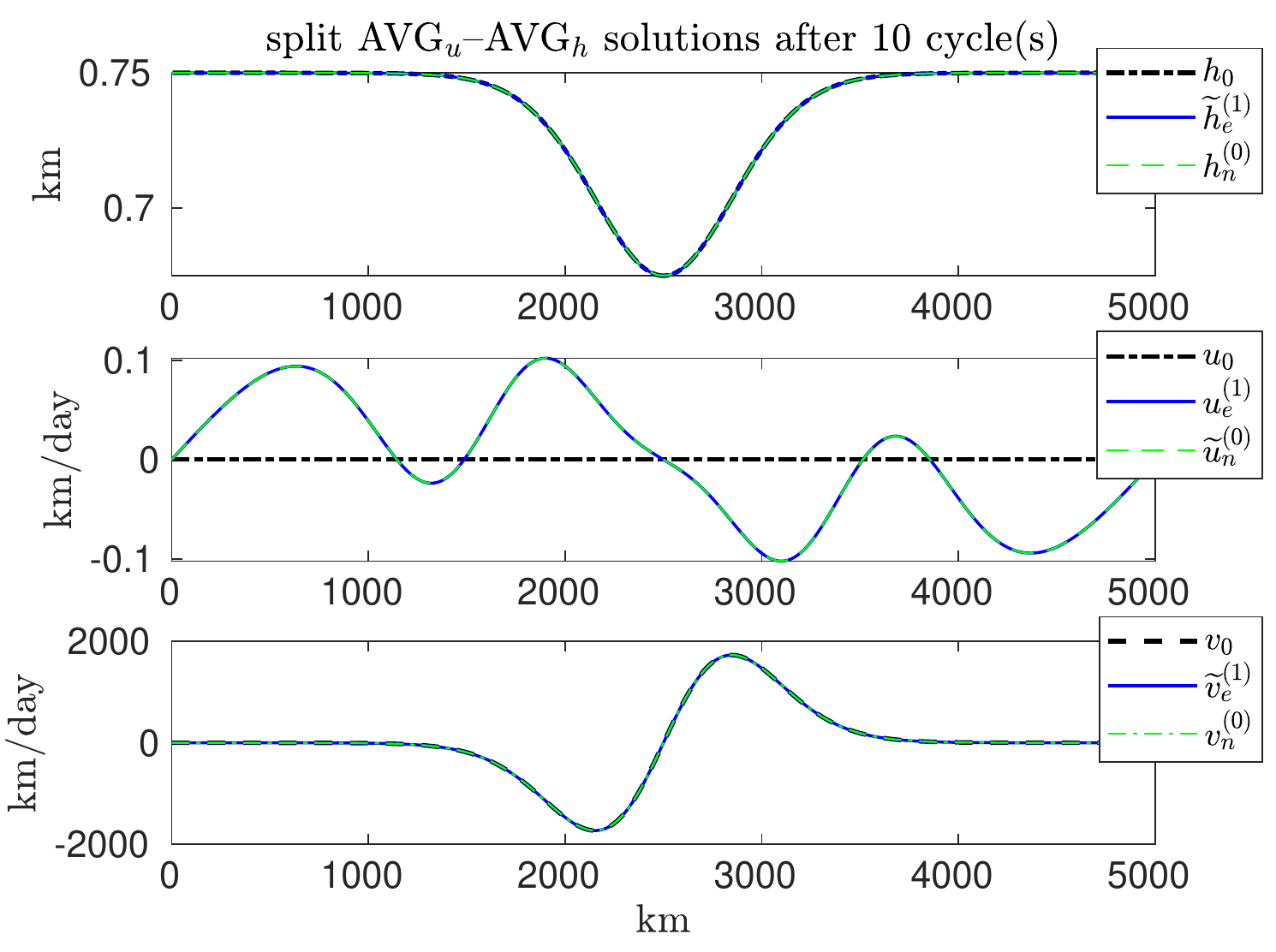}}  & 
    \hspace{0cm}{\includegraphics[scale=0.37]{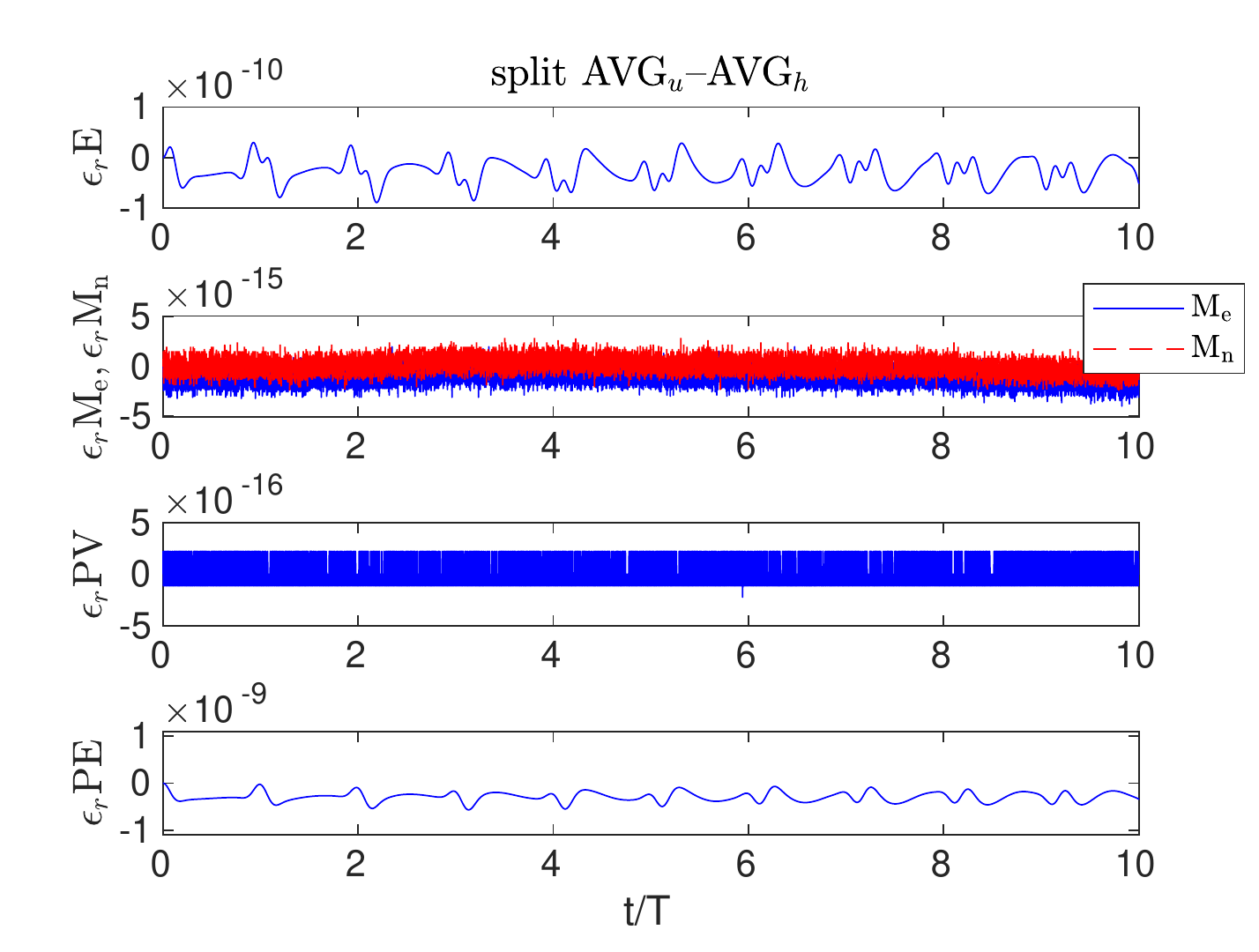}}  
   \end{tabular}
   \vspace{-0.5em}
    \caption{Solutions for the split $\AVGu$-$\,\AVGh$ scheme for a mesh with 512 elements. Initial fields are shown as dashed-dotted lines.
    Left: flow in geostrophic balance after 10 cycles. 
    Right: time series of the quantities of interest for 10 cycles. 
    }
    \label{fig_tc1_H075_H10fields}
    \vspace{-2em}
  \end{figure}

The study of structure preservation (topological properties) will be performed with a flow
in geostrophic balance in which the terms are linearly balanced 
while nonlinear effects are comparably small (Fig.~\ref{fig_tc1_H075_H10fields}). 
To illustrate the long term behaviour, we run the simulation in this TC 1 for about 10 cycles 
(meaning that the (analytical) wave solutions have traveled 10 times over 
the entire domain). 
To test convergence and stability (i.e. metric-dependent properties),
we use in TC 2 a steady state solutions of Eqn.~\eqref{slice-model}. 
To illustrate the metric dependency of the dispersion relations, 
we use in TC 3 an initial height distribution (as in Fig.~\ref{fig_tc1_H075_H10fields}, left) 
that is only partly in linear geostrophic balance such that shock waves 
with small scale oscillations develop that depend on the dispersion relation. 
More details on the TC can be found in \cite{BaBeCo2019}.

\vspace{-0.5em}

\begin{figure}[h] \centering
  \begin{tabular}{cc} 
   \hspace{0cm}\includegraphics[scale=0.33]{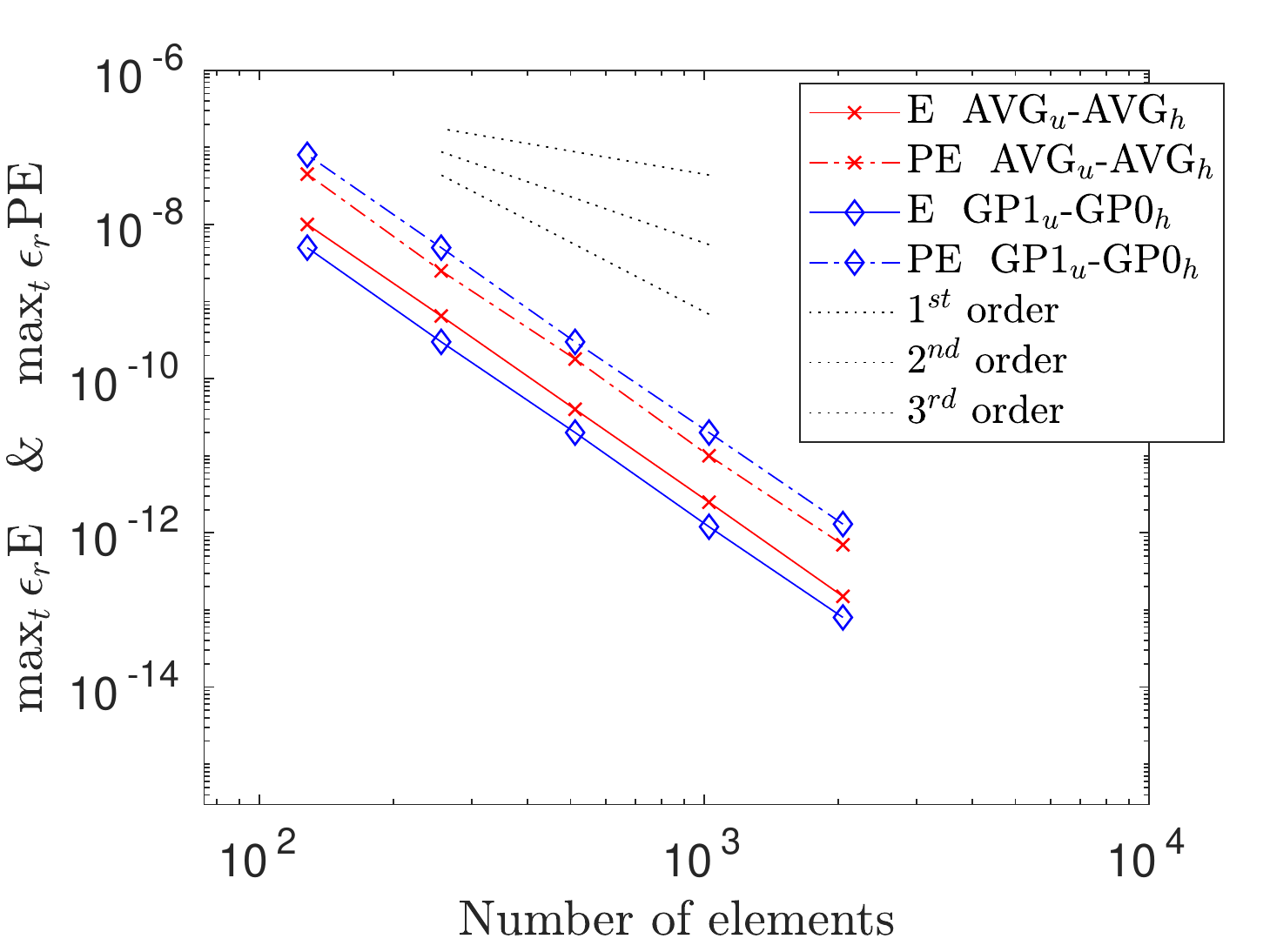}   &
    \hspace{0.5cm}\includegraphics[scale=0.33]{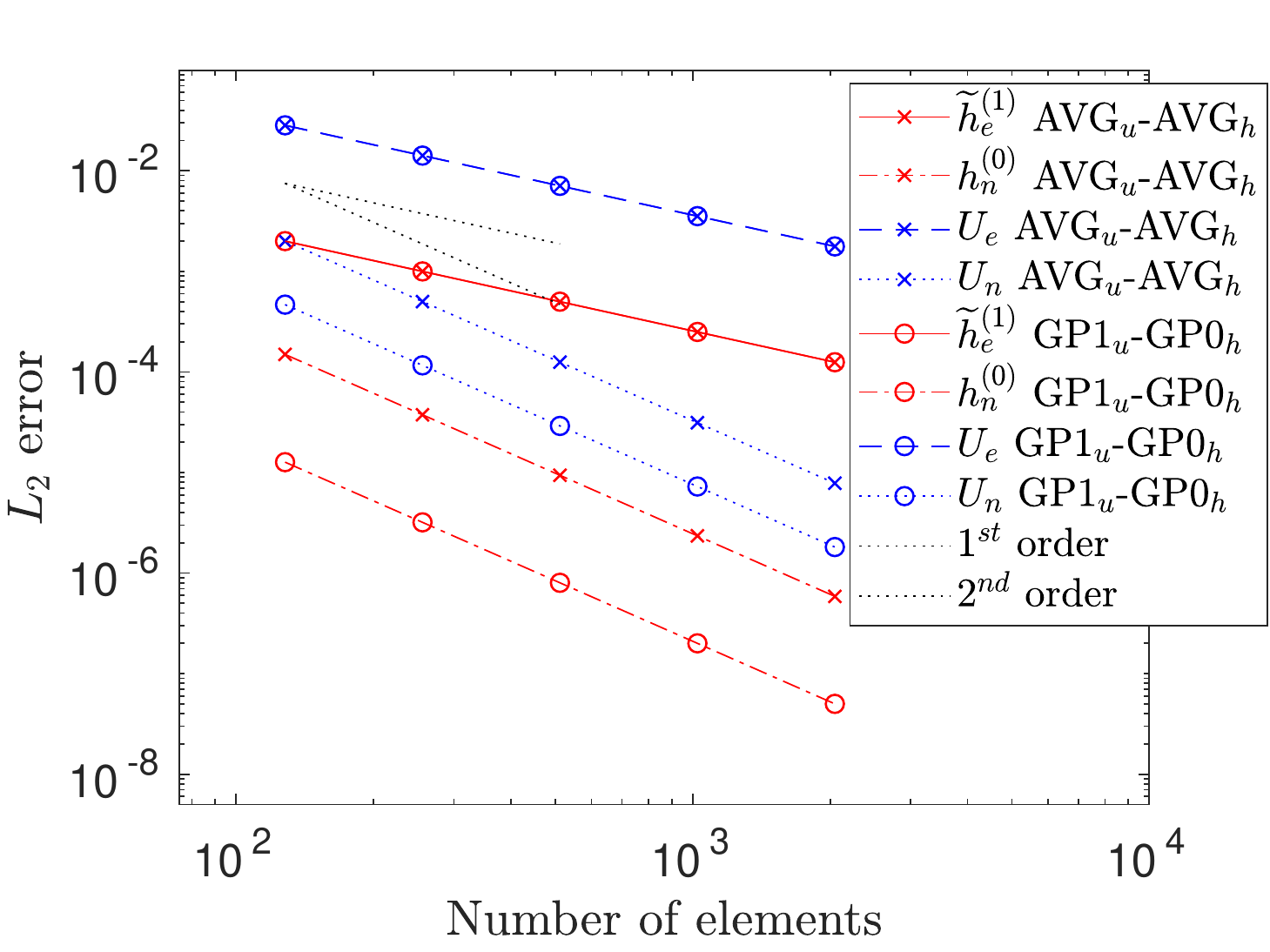} 
  \end{tabular}
  \vspace{-0.5em}
  \caption{Relative error values in dependence of $N$ of $\AVGu$-$\,\AVGh$ compared to $\GPIu$--$\,\GPOh$.
  Left: errors for $E$ and $PE$ for TC 1. 
  Right: errors for the steady state solution of TC 2 after 1 cycle.}
\label{fig_tc1_diag}
\vspace{-2em}
  \end{figure}

\smallskip\noindent\textbf{Topological properties.}
Figure~\ref{fig_tc1_diag} (left) shows for TC 1 the relative errors of the averaged 
split scheme for energy $E$, mass $M_e$ or $M_n$, potential vorticity $PV$ and 
enstrophy $PE$ (see definitions in \cite{BaBeCo2019}).
In all cases studied, $M_e, M_n$ and $PV$ are preserved at machine precision. 
For $512$ elements, the error in $E$ is at the order of $10^{-10}$ while
for $PE$ it is at $10^{-9}$. With higher mesh resolution, the errors decrease with 
third order rate (Fig.~\ref{fig_tc1_diag}, left).

When compared to the split schemes of \cite{BaBeCo2019}, these error values are very 
close to the results presented therein, underpinning the fact
that modifications in the metric equations do not affect the quality of 
structure preservation of the schemes. 
 
 \vspace{-1.5em}

\smallskip\noindent\textbf{Metric-dependent properties.}
Consider next the convergence behaviour of the averaged split scheme $\AVGu$-$\,\AVGh$
shown in Figure~\ref{fig_tc1_diag} (right) for TC 2. To ease comparison, 
we include $L_2$ error values of the split scheme $\GPIu$--$\,\GPOh$
of \cite{BaBeCo2019} noting that the other split schemes presented therein 
share more or less the same error values for the corresponding fields.
In all cases, the error values decrease as expected: all P1 fields 
show second order, all P0 fields first order convergences rates. 

While the errors of the P0 fields of $\AVGu$--$\,\AVGh$
is close to the corresponding values of the split schemes of \cite{BaBeCo2019},
the P1 fields of $\AVGu$--$\,\AVGh$
have error values that are about one order of magnitude large
than the corresponding fields of e.g. $\GPIu$--$\,\GPOh$.
This agrees well with the fact that we do not solve the full linear system 
in the metric equations to recover the P1 fields but use instead approximations, 
which slightly increases the P1 error values of $\AVGu$--$\,\AVGh$. 

\smallskip
 \begin{figure}[h!] \centering
  \begin{tabular}{ccc} 
  \hspace*{0cm}\begin{overpic}[scale=0.3,unit=1mm]
    {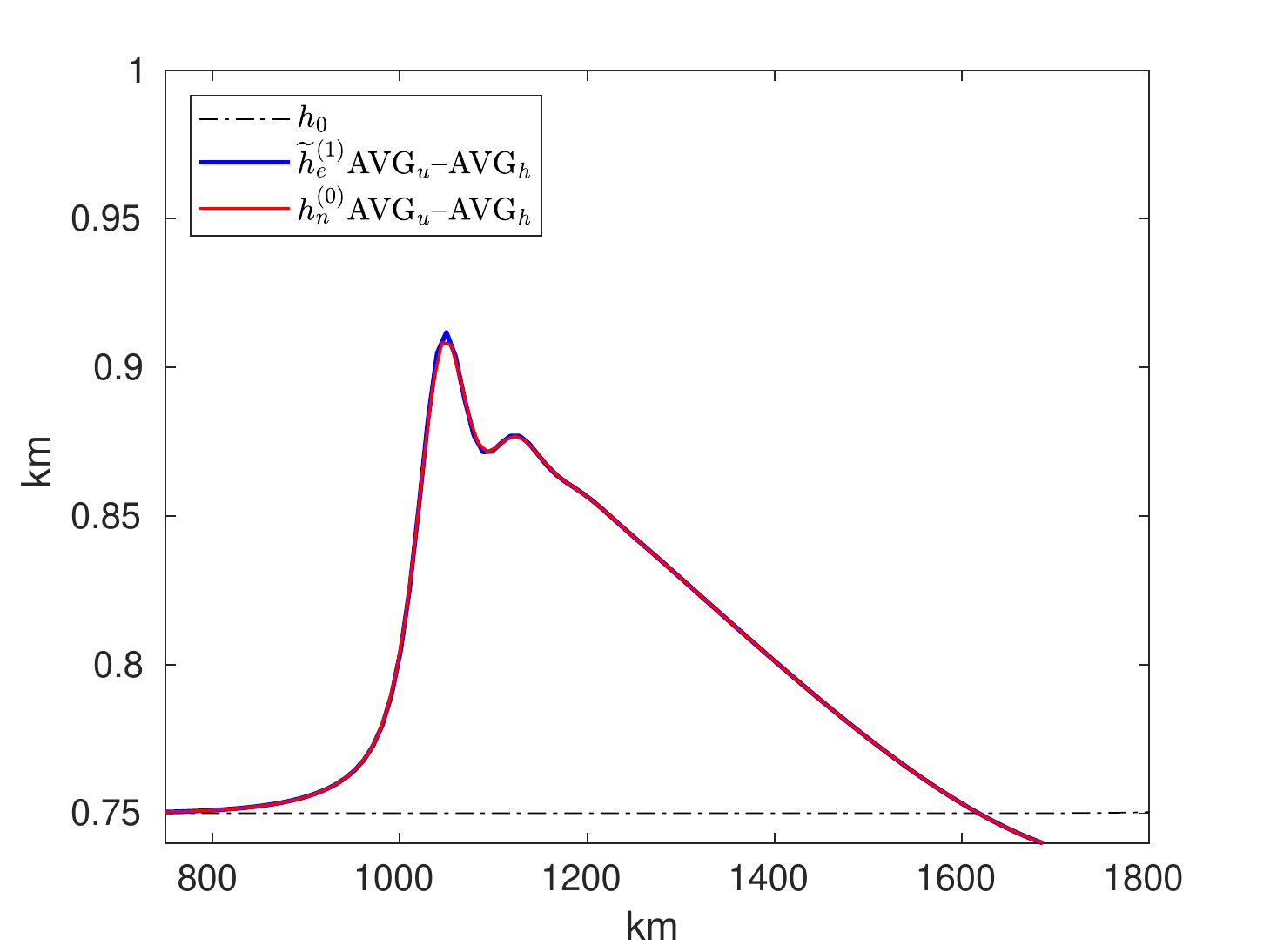}
     \put(45,40){\includegraphics[scale=.205]{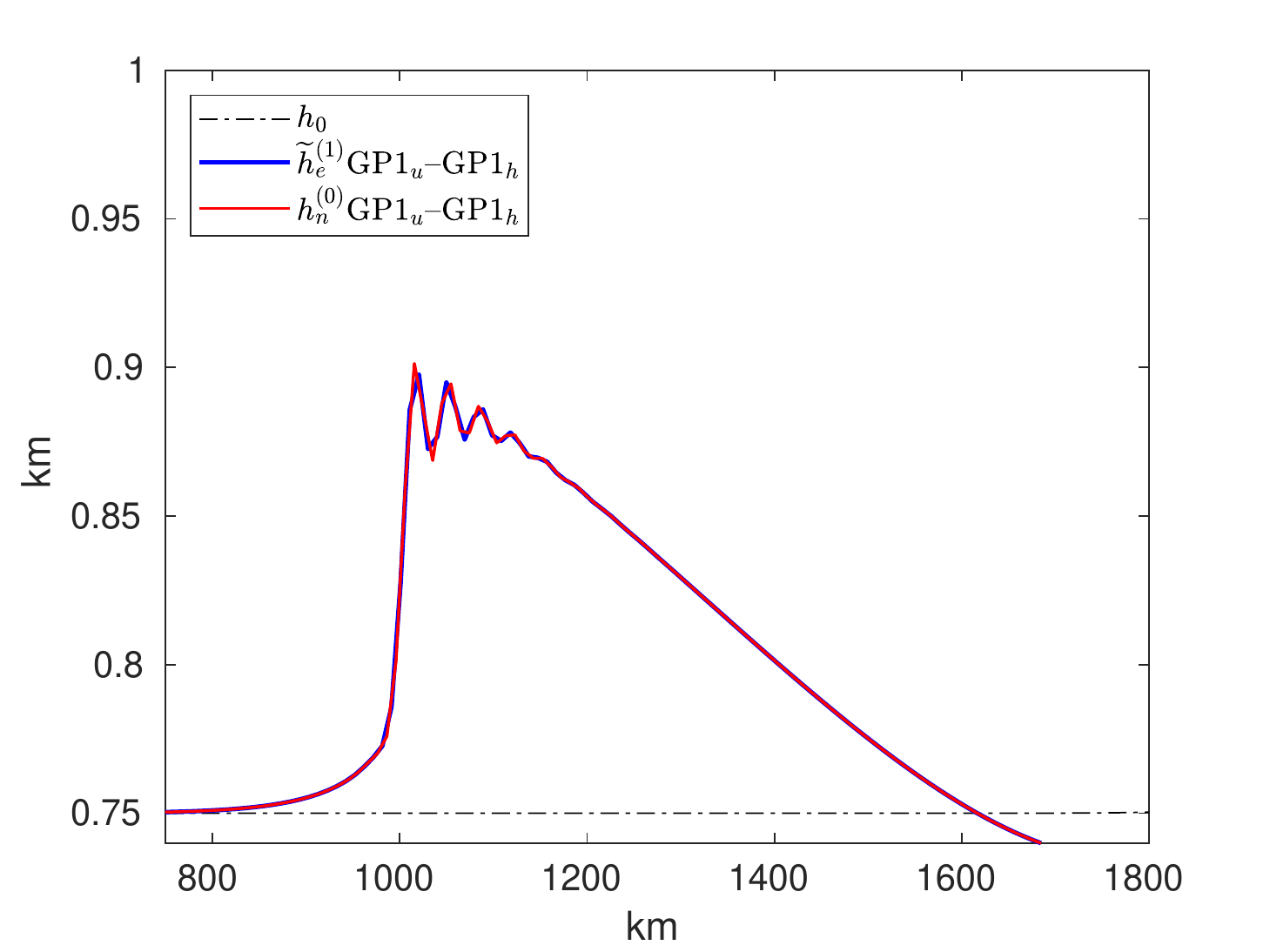}}
     \put(80,20){\includegraphics[scale=.205] {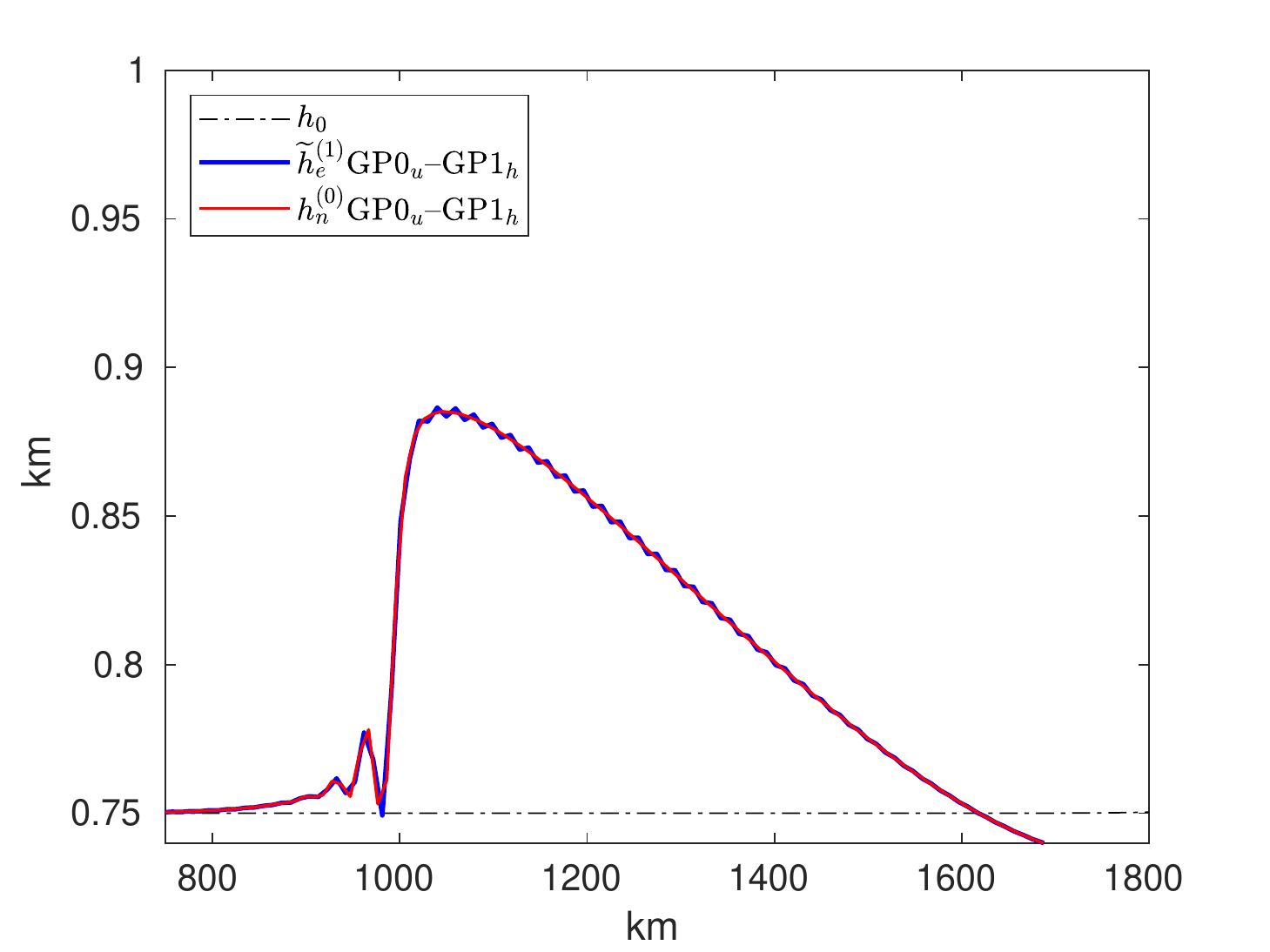}}
     \put(115,0){\includegraphics[scale=.205] {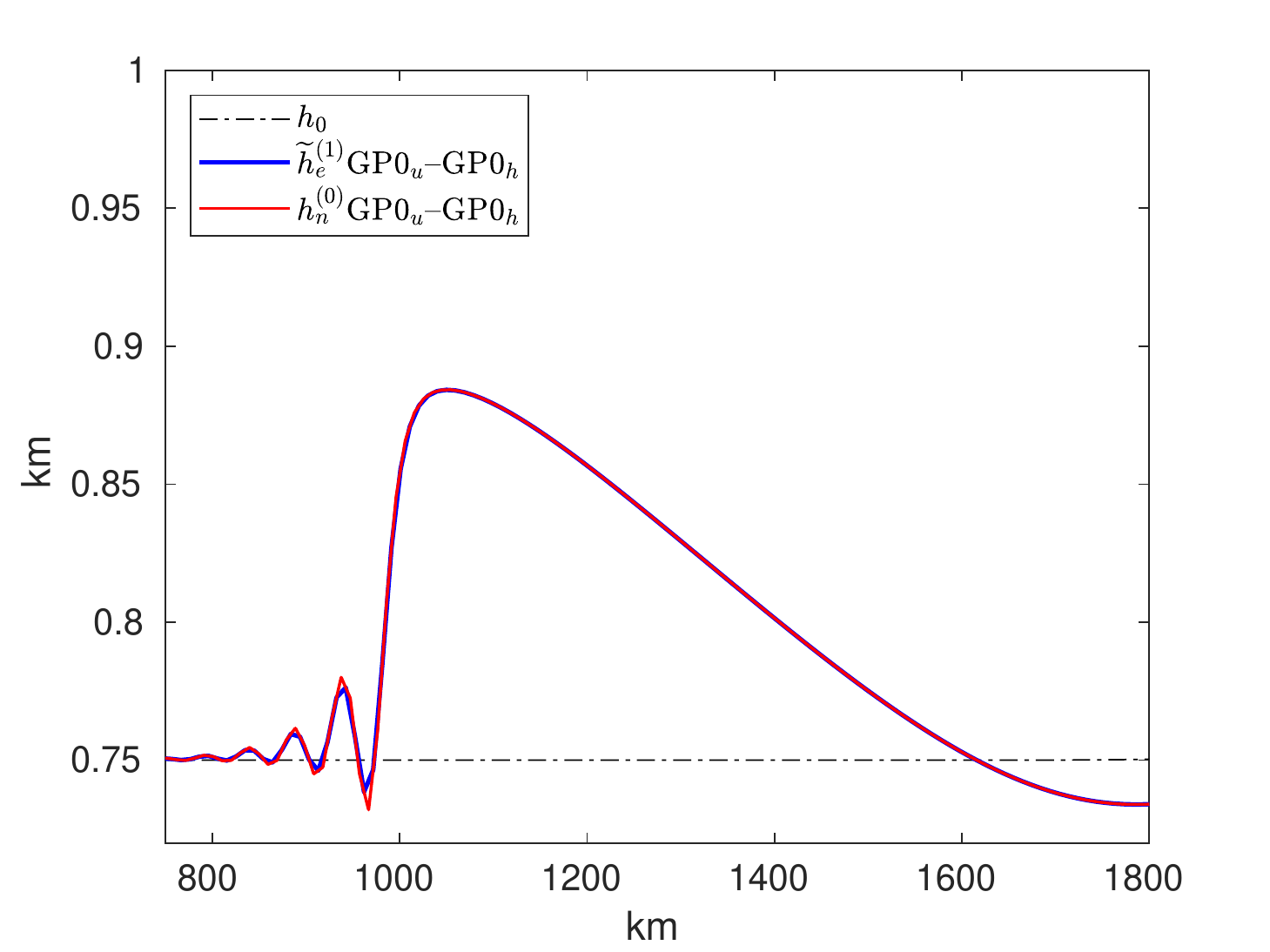}}
     \put(15,50){${\color{magenta} \omega_{aa}}$}
     \put(155,40){${\color{red} \omega_{00}}$}
     \put(125,60){${\color{blue} \omega_{01}}$}
     \put(90,80){${\color{green} \omega_{11}}$}
   \end{overpic}
   \hspace{2.8cm} \begin{overpic}[scale=0.30,unit=1mm]
   {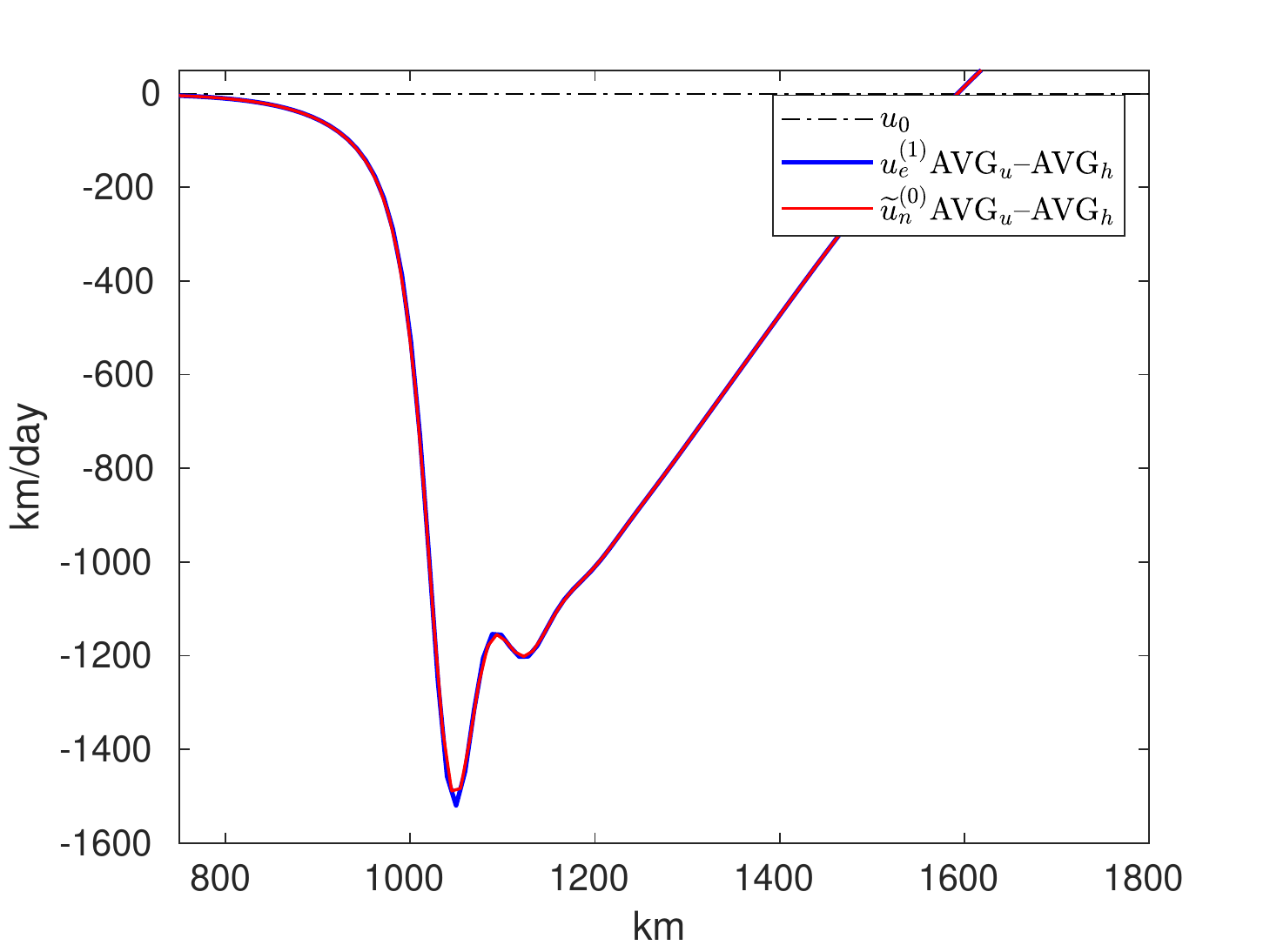}
   \end{overpic}
   \end{tabular}
   \vspace{-0.5em}
   \caption{Fields with oscillations at the wave fronts in dependency of the 
   wave dispersion relations of Fig.~\ref{fig_disprel} on a mesh with $N_e = 512$ 
   elements and after a simulation time of $0.225$ cycles.
   }  
    \label{fig_tc3_disprel1}
    \vspace{-2em}
  \end{figure}

 With TC 3 we illustrate numerically how the choice of metric equations determines the 
 discrete dispersion relations.  
 As derived in Sect.~\ref{sec_numerics}, 
 the discrete dispersion relation of $\AVGu$--$\,\AVGh$ equals a sine wave,
 hence all waves of frequency $k$ have wave speeds equal or slower than $c$ 
 (black curve in Fig.~\ref{fig_disprel}). In particular for wave numbers 
 larger then $\frac{\pi}{2 \Delta x}$, waves start to slow down until
 there is a standing wave at $k = \frac{\pi}{\Delta x}$.
 This is a similar behavior to the $\GPIu$--$\,\GPIh$ scheme of \cite{BaBeCo2019}, 
 but for $\AVGu$--$\,\AVGh$ this effect is stronger given the generally slower 
 wave propagation. This behaviour is clearly visible in Figure~\ref{fig_tc3_disprel1}
 where we observe in both fields lower frequency oscillation behind the 
 front when compared to $\GPIu$--$\,\GPIh$ (see inlet). This result agrees 
 well with the discrete dispersion relations shown in Fig.~\ref{fig_disprel}.

 \vspace{-2em}
 \section{Conclusions}
 \label{sec_conclusion}
 \vspace{-1em}
 
 We introduced a y-independent RSW slice-model in 
 split Hamiltonian form and derived a family of 
 lowest-order (P0-P1) structure-preserving split schemes.
 The splitting of the equations into topological and metric parts 
 transfers also to schemes' properties. The framework allows 
 for different realizations of metric equations which all 
 preserve the Hamiltonian and the Casimirs of the Poisson bracket. 
 This allowed us to introduce an approximation of the metric equations 
 which is structure-preserving, achieving a speedup of 
 a factor of 2 because no linear systems had to be solved.

\end{document}